\documentclass[12pt,a4paper,twoside]{amsart}

\usepackage[T1]{fontenc}
\input xypic\usepackage{amsmath}
\usepackage{amssymb}

\usepackage{latexsym}

\usepackage{graphicx}
\usepackage{psfrag}
\usepackage{psfig}

\begin{document}
\newcommand{\refeq}[1]{(\ref{#1})}

\def\tilde{\widetilde}\newcommand{\HOX}[1]{\marginpar{\footnotesize #1}}
\def\hat{\widehat}\def\cal{\mathcal}
\def \beq {\begin {eqnarray}}
\def \eeq {\end {eqnarray}}
\def \ba {\begin {eqnarray*}}
\def \ea {\end {eqnarray*}}
\def \bfo {\ba}
\def \efo {\ea}
\def\Om{\Omega}
\def\O{{\mathcal O}}

\newtheorem{theorem}{Theorem}[section]
\newtheorem{lemma}[theorem]{Lemma}
\newtheorem{Mtheorem}{Theorem}
\newtheorem{proposition}[theorem]{Proposition}
\newtheorem {adefinition}[theorem]{Assumption}
\newtheorem {definition}[theorem]{Definition}
\newtheorem {example}[theorem]{Example}
\newtheorem {corollary}[theorem]{Corollary}
\newtheorem {problem}[theorem]{Problem}
\newtheorem {ex}[theorem]{Exercise}

\def\proof{\noindent{\bf Proof.}\quad}

   \def \cA {{\cal A}}\def\vrho{\varrho}
\def \fg{{\bf g}}
\def \cB {{\cal B}}
\def \cD {{\cal D}}
\def \A {{\cal A}}
\def \B {{\cal B}}
\def \D {{\cal D}}
\def \Disc {{\Bbb D}}
\def \tX {{\widetilde X}}
\def \t {{\tau}}
\def \Box {\quad $\square$\medskip }
\def \ssskip {\vskip -23pt \hskip 40pt {\bf.} }
\def \sskip {\vskip -23pt \hskip 33pt {\bf.} }
\def\a{\alpha}
\def \e{\varepsilon}
\def \b{\beta}
\def \bn{\underline n}
\def\Grad{\hbox{Grad}\,}\def\Div{\hbox{Div}\,}
\def \cW {{\cal W}}
\def\Hess{{\rm Hess}\,}
\def\sign{{\rm sign}\,}
\def\re{\hbox{Re}\,}
\def\im{\hbox{Im}\,}
\def\Re{\hbox{Re}\,}         
\def\Im{\hbox{Im}\,}

\def\proof{\noindent{\bf Proof.}\quad}

\title{CALDER\'ON'S INVERSE PROBLEM FOR ANISOTROPIC CONDUCTIVITY IN THE PLANE}
\author {Kari Astala}
\address{Rolf Nevanlinna Institute, University of Helsinki,
P.O.~Box~4 (Yliopistonkatu~5), FIN-00014 University of Helsinki,
Finland}
\email{Kari.Astala@helsinki.fi, Matti.Lassas@helsinki.fi,
$ \ \ \ \ \ \ \ \ \ \ \ \ \ \ \ \ \ \ \ \ $  \linebreak
\hbox{ } $ \ \ \ \ \ \ \  \ \ \ \ \ \ \ $ljp@rni.helsinki.fi}

\author{Matti Lassas}

\author{Lassi P\"aiv\"arinta}
\maketitle

\def\Hs{H^{s+1}_0(\partial M\times [0,T])}
\def\H2s{H^{s+1}_0(\partial M\times [0,T/2])}
\def\smooth {C^\infty_0(\partial M\times [0,2r])}
\def\supp{\hbox{supp }}
\def\diam{\hbox{diam }}
\def\dist{\hbox{dist}}
\def\R{{\mathbb R}}
\def \fR{\R}
\def\Z{{ \mathbb Z}}
\def\C{{\mathbb  C}}
\def\e{\varepsilon}
\def\g{\tau}
\def\F{{\cal F}}
\def\N{{\cal N}}
\def\U{{\cal U}}
\def\W{{\cal W}}
\def\O{{\cal O}}
\def\exp{\text{exp}}

\def\l{\sigma}
\def\o{\over}
\def\i{\infty}
\def\d{\delta}
\def\G{\Gamma}
\def\exp{\hbox{exp}}
   \def\bra{\langle}
\def\cet{\rangle}
\def\p{\partial}

\def\tM{\tilde M}
\def\tR{\tilde R}
\def\tr{\tilde r}
\def\tz{\tilde z}
\def\expec{\hskip1pt{\mathbb  E}\hskip2pt}
\def\prob{{\cal P}}
\def\good{{\cal A}}
\def\proof{\noindent{\bf Proof.}\quad }
\def\ds{\displaystyle}
\def\pz{\p_z}
\def\pbz{\p_{\overline  z}}
\def\bar{\overline}
\def\H{{\cal H}}

{\bf Abstract:} We study inverse conductivity problem for an
anisotropic conductivity  $\sigma\in  L^\infty$ in  bounded and unbounded domains.
Also, we give applications of the results in the case when 
Dirichlet-to-Neumann and
Neumann-to-Dirichlet maps are given only on a part of the boundary.

\section{INTRODUCTION}
Let us consider the 
anisotropic conductivity equation in two dimensions
\beq\label{conduct}
\nabla\cdotp \sigma\nabla u= \sum_{j,k=1}^2 
\frac \p{\p x^j}\sigma^{jk}(x) \frac \p{\p x^k} u &=& 0\hbox{ in } \Omega,\\
u|_{\p \Omega}&=&\phi.\nonumber
\eeq
Here 
 $\Omega\subset\R^2$ is a simply connected domain. 
The
 conductivity
$\sigma=[\sigma^{jk}]_{j,k=1}^2$ is  a
symmetric, positive definite matrix function, and
$\phi \in H^{1/2}(\p \Omega)$ is the prescribed voltage on the boundary.
Then it is well known that equation (\ref{conduct}) has
a unique solution $u\in H^1(\Om)$.

In the case when $\sigma$ and $\p \Omega$
are smooth, we can define the voltage-to-current
(or Dirichlet-to-Neumann) map by 
\beq
\Lambda_\sigma(\phi)= Bu|_{\p \Omega}
\eeq
where
\beq
Bu=\nu \cdotp \sigma \nabla u,
\eeq
 $u\in H^1(\Om) $ is the solution of (\ref{conduct}), and $\nu$
is the unit 
normal vector of $\p\Omega$.
Applying
the divergence theorem, we have
\beq\label{Q_l}
Q_{\l,\Om} (\phi):=\int_\Omega  \sum_{j,k=1}^2\sigma^{jk}(x)\frac {\p u}{\p x^j}
\frac {\p u}{\p x^k}
dx=\int_{\p \Omega} \Lambda_\l(\phi) \phi\, dS,
\eeq
where $dS$ denotes the arc lenght on $\p \Omega$.
The quantity $Q_{\l,\Om} (\phi)$ represents the
power needed to maintain the potential $\phi$ on $\p \Omega$. By symmetry
of $\Lambda_\sigma$,
knowing $Q_{\l,\Om}$ is
equivalent with knowing
$\Lambda_\l$. For general $\Omega$ and $\sigma\in L^\infty(\Om)$,
the trace $u|_{\p \Om}$ is defined as the equivalence class
of $u$ in $H^1(\Om)/H^1_0(\Om)$ (see \cite{AP}) and
formula (\ref{Q_l}) is used to define the map $\Lambda_\sigma$.

If $F:\Omega\to \Omega,\quad F(x)=(F^1(x),F^2(x))$, is a
diffeomorphism with $F|_{\p
\Omega}=\hbox{Identity}$, then by making the change of variables $y=F(x)$
and setting $v=u\circ F^{-1}$ in the first
integral in (\ref{Q_l}), we obtain
\ba
\nabla\cdotp (F_*\sigma)\nabla v=0\quad\hbox{in }\Om,
\ea
where
\beq\label{cond and metr}
(F_*\sigma)^{jk}(y)=\left.
\frac 1{\det [\frac {\p F^j}{\p x^k}(x)]}
\sum_{p,q=1}^2 \frac {\p F^j}{\p x^p}(x)
\,\frac {\p F^k}{\p x^q}(x)  \sigma^{pq}(x)\right|_{x=F^{-1}(y)},
\eeq
or
\beq\label{5 1/2}
F_*\sigma(y)=\left.\frac 1{J_F(x)} DF(x)\,\sigma(x)\,DF(x)^t\right|_{x=F^{-1}(y)},
\eeq
is the push-forward of the conductivity $\l$ by $F$.
Moreover, since $F$ is identity at $\p \Om$, we obtain
from (\ref{Q_l}) that
\ba
\Lambda_{F_*\l}=\Lambda_\l.
\ea
Thus, the change of coordinates shows that there is a large
 class of conductivities which give rise to the
same electrical
measurements at the boundary.

We consider here the converse question, that if we have two
conductivities which have the same Dirichlet-to-Neumann map, is it the case 
that each of them can be obtained by 
pushing forward the other.

In applied terms,
this inverse problem to determine $\sigma$ (or its properties)
from $\Lambda_{\sigma}$ is
 also known as {\em Electrical Impedance Tomography}. 
It has been proposed as a valuable diagnostic, see \cite{CIN99}.

 In the case where $\sigma^{jk}(x)=\sigma(x)\delta^{jk}$, $\sigma(x)\in \R_+$,
the metric is said to be isotropic. 
In 1980 it was proposed by A.~Calder\'on
 \cite {Cl} that in the isotropic case any
bounded conductivity $\sigma(x)$ might be determined
solely from the boundary measurements, i.e., from $\Lambda_\sigma$.
Recently this has been confirmed in the two dimensional case
(c.f. \cite{AP}). In the case when isotropic $\sigma$ is smoother
than just a $L^\infty$-function, the same conclusion
is known to hold also in higher dimensions. 

The first global
uniqueness result was obtained for a $C^\infty$--smooth
conductivity in dimension $n\geq 3$ by 
 J.~Sylvester and G.~Uhlmann  in 1987 \cite{SylUhl}. 
In dimension two A.~Nachman \cite{Nach2D} produced in 1995 a uniqueness
result for conductivities with two derivatives. 
The corresponding algorithm has been successfully implemented and 
proven to work efficiently even with real data \cite{Siltanen1,Siltanen2}. 
The  reduction of 
regularity assumptions has since been under active study.
In dimension two the optimal $L^\infty$-regularity
was obtained in \cite{AP}.
In dimension $n\geq 3$ the uniqueness has presently been shown for
isotropic conductivities $\sigma \in W^{3/2,\infty}(\Om)$
in \cite{PPU}
and for globally
 $C^{1+\e}$--smooth isotropic conductivities  having only
co-normal singularities in \cite{GLU1}. 

Also, the stability of reconstructions of the inverse conductivity problem
have been extensively studied. For these results, see \cite{Al,AlS,BBR}
where stability results are based on reconstruction
techniques of \cite{BU} in dimension two and
those of  \cite{Na1} in dimensions $n\geq 3$. 

In anisotropic case, where $\sigma$ is a  matrix function and the
problem is to recover the conductivity $\sigma$ up to the action
of a class of diffeomorphisms, much
less is known. In dimensions $n\geq 3$ it is generally known  
only that piecewise analytic conductivities can be constructed
(see \cite{KV1,KV2}). For Riemannian manifolds this kind of
technique has been generalized in \cite{LeU,LaU,LTU}.
In dimension $n=2$ the inverse problem has  been considered 
by J.~Sylvester   \cite{Sy1} for $C^3$ and Z.~Sun and G.~Uhlmann \cite{SU}
for $W^{1,p}$-conductivities. The idea of 
\cite{Sy1} and \cite{SU} is that under quasiconformal change
of coordinates (cf. \cite{Alf,IM}) any anisotropic
conductivity can be changed to isotropic one, see also
section \ref{Sec: proofs} below. The 
purpose of this paper is to carry this technique over
to the $L^\infty$--smooth case and then use the result of
\cite{AP} to obtain uniqueness up to the group of diffeomorphisms.

The advantage of the reduction of the smoothness assumptions
up to $L^\infty$
does not lie solely on the fact that many  conductivities
have jump-type singularities but it also allows us to consider 
much more complicated singular structures such as porous rocks
\cite{Che}. Moreover it is important that this approach 
 enables us to consider general diffeomorphisms.
Thus anisotropic inverse problems in half-space or exterior 
domains can be solved simultaneously. This will be considered
in Section \ref{sec: con}.

If $\Om\subset \R^2$ is a bounded domain, it is convenient
to consider the class of  matrix functions $\sigma=[\sigma^{jk}]$
such that 
\beq\label{basic ass}
[\sigma^{ij}]\in L^\infty(\Om;\R^{2\times 2}),
\quad [\sigma^{ij}]^t=[\sigma^{ij}],\quad
C_0^{-1}I \leq [\sigma^{ij}]\leq C_0I
\eeq
where $C_0>0$. In sequel, the minimal possible value of $C_0$
is denoted by $C_0(\sigma)$.
We use the notation
\ba
\Sigma(\Om)=\{\sigma\in L^{\infty}(\Om;\R^{2\times 2})& |&
\ C_0(\sigma)<\infty\}.
\ea
Note that it is necessary
to require $C_0(\sigma)<\infty$ 
as otherwise there would be counterexamples showing
that even the equivalence class of the conductivity can not be recovered
 \cite{GLU2a,GLU2}.

Our main goal in this paper is to show that
an anisotropic $L^\infty$--conductivity can be determined
up to a $W^{1,2}$-diffeomorphism:

\begin{Mtheorem}\label{theorem1}
Let $\Om\subset\R^{2}$ be a simply connected bounded domain and 
$\sigma\in L^{\infty}(\Om;\R^{2\times 2})$. Suppose that
the assumptions (\ref{basic ass}) are valid. Then
the Dirichlet-to-Neumann map 
$\Lambda_{\sigma}$
determines the equivalence class
\ba
E_\sigma=\{\sigma_1\in \Sigma(\Om)& |& 
\hbox{$\sigma_1 = F_*\sigma$, 
$F:\Om\to \Om$ is  $W^{1,2}$-diffeomorphism and}\\
& &F|_{\p \Om}=I\}.
\ea
\end{Mtheorem}

We prove this result in Section \ref{Sec: proofs}.

Finally, note that the $W^{1,2}$--diffeomorphisms $F$
preserving the class $\Sigma(\Om)$ are precisely the quasiconformal
mappings. Namely, if $\sigma_0\in \Sigma(\Om)$ and $\sigma_1=F_*(
\sigma_0)\in \Sigma(F(\Om))$ then
\beq\label{Kari 8}
\frac 1{C_0}||DF(x)||^2I\leq DF(x)\, \sigma_0(x)\, DF(x)^t\leq 
C_1J_F(x)I
\eeq
where $I=[\delta^{ij}]$ and we obtain
\beq\label{Kari 9}
||DF(x)||^2\leq K J_F(x),\quad \hbox{for a.e. }\ x\in \Om
\eeq
where $K=C_1C_0<\infty$. Conversely, if (\ref{Kari 9}) holds
and $F$ is $W^{1,2}_{loc}$-homeomorphism then
$F_*\sigma\in \Sigma(F(\Om))$ whenever $\sigma\in \Sigma(\Om)$.
Furhtermore, recall that a map $F:\Om \to \tilde \Om$ is quasiregular
if $F\in W^{1,2}_{loc}(\Om)$ and the condition (\ref{Kari 9}) holds.
Moreover, a map $F$ is quasiconformal if it is
quasiregular and a $W^{1,2}$--homeomorphism.

\section{CONSEQUENCES AND APPLICATIONS OF THEOREM \ref{theorem1}}\label{sec: con}

Here we consider applications of the diffeomorphism-technique
to various inverse problem. The formulated results,
Theorems \ref{Lem: A1}--\ref{Lem: A3} are proven in
Section \ref{sec: proof of con}.

\subsection{Inverse Problem in the Half Space}

Inverse problem in half space is of crucial importance
in geophysical prospecting, 
seismological imaging, non-destructive testing etc.
For instance, the imaging of soil was the original
motivation of Calder\'on's seminal paper \cite {Cl}.
As we can use a diffeomorphism to map the open half space
to the  unit disc, we can apply the previous result for the half space case.
One should observe that in this deformation even
infinitely smooth conductivities can become  
non-smooth at the boundary 
(e.g.  conductivity oscillating near infinity produces
a non-Lipschitz conductivity in push-forward) and thus
the  low-regularity result \cite{AP} is essential for the problem.

Thus, for $\sigma\in \Sigma(\R^2_-)$ let us consider the problem
\beq\label{conduct D1a}
\nabla\cdotp \sigma\nabla u&=& 0\hbox{ in } \R^2_-=\{(x^1,x^2)\,|\, x^2<0\}
,\\
u|_{\p \R^2_-}&=&\phi,\label{conduct D1b}\\
u&\in& L^\infty(\R^2_-)\label{conduct D1c}.
\eeq
Notice that here  the radiation condition at infinity (\ref{conduct D1c})
is quite simple. We assume just that the potential $u$ 
does not blow up at infinity.
The equation (\ref{conduct D1a}--\ref{conduct D1c}) is uniquely solvable 
and as before we can define
\ba
\Lambda_\sigma:
H_{comp}^{1/2}(\p \R^2_-)\to H^{-1/2}(\p \R^2_-),\quad
\phi\mapsto \nu\cdotp\sigma\nabla u|_{\p \R^2_-}.
\ea
\begin{theorem}\label{Lem: A1} The map $\Lambda_\sigma$
 determines the 
equivalence class
\ba
E_{\sigma}=\{\sigma_1\in \Sigma(\R^{2}_-)& |& 
\hbox{$\sigma_1 = F_*\sigma$, 
$F:\R^2_-\to \R^2_-$ is  $W^{1,2}$-diffeomorphism,}\\
& &F|_{\p \R^2_-}=I\}.
\ea
Moreover, each orbit $E_{\sigma}$ contains at most
one isotropic conductivity, and consequently 
if $\sigma$ is known to be isotropic, it is determined uniquely by $\Lambda_
\sigma$.
\end{theorem}
Note that the natural growth requirement 
$ \lim_{|z|\to \infty} |F(z)|=\infty$ 
follows automatically from the above assumptions on $F$.

\subsection{Inverse Problem in the Exterior Domain}

An inverse problem similar to that of the half space can be considered
in an exterior domain where one wants to find the conductivity in
a complement of
a bounded simply connected domain. This type of problem is encountered in
cases where measurement devices are embedded to an unknown domain.

In the case of $S=\R^2\setminus \overline D$,
 where $D$ is a bounded Jordan domain,
we consider the problem
\beq\label{conduct S}
\nabla\cdotp \sigma\nabla u&=& 0\quad \hbox{ in } S,\\
u|_{\p S}&=&\phi\in H^{1/2}(\p S),\label{conduct S1} \\
u&\in& L^\infty(S).\label{conduct S2a}
\eeq
Again, the radiation condition (\ref{conduct S2a}) of infinity
is only that the solution is uniformly bounded.
For this equation we define
\ba
\Lambda_\sigma:
H^{1/2}(\p S)\to H^{-1/2}(\p S),\quad
\phi\mapsto \nu\cdotp\sigma \nabla u|_{\p S}.
\ea
Surprisingly, the result is different from the half-space case.
The reason for this is the phenomenon that 
the group of diffeomorphisms preserving the data does not fix the point
of the infinity. More precisely,
there are two points $x_0,x_1\in S\cup\{\infty\}$  
such that $F(x_0)=\infty$, $F^{-1}(x_1)=\infty$, and
$F:S\setminus \{x_0\}\to S\setminus \{x_1\}$. 
In particular, this means that the uniqueness does not hold up to
diffeomorphisms mapping the exterior domain to itself.

For convenience, we compactify $S$ by adding one infinity point,
denote $\overline S=S\cup \{\infty\}$, and define $\sigma(\infty)=1$.
We say that $F:\overline S\to \overline S$ is a $W^{1,2}$-diffeomorphism
if $F$ is homeomorphism and a $W^{1,2}$-diffeomorphism
in spherical metric \cite{Alf}.

\begin{theorem}\label{Lem: A2} Let $\sigma\in \Sigma(S)$. 
Then the map $\Lambda_\sigma$ 
 determines the 
equivalence class
\ba
E_{\sigma,S}=\{\sigma_1\in \Sigma(S)& |& 
\hbox{$\sigma_1 = F_*\sigma$, 
$F:\overline S\to \overline S$ is a  $W^{1,2}$-diffeomorphism,}\\
& &F|_{\p S}=I\,\}.
\ea
Moreover, if $\sigma$ is known to be isotropic, it is determined uniquely
by  $\Lambda_\sigma$.
\end{theorem}

\subsection{Data on Part of the Boundary}

In many inverse problems data is measured
only on a part of the boundary. For the conductivity
equation in dimensions $n\geq 3$ it
has been shown that if the measurements are done
on a part of the boundary, then the integrals of the unknown conductivity over
certain 2-planes can be determined \cite{GU}.
In one-dimensional inverse problems partial data is 
often considered with two different boundary conditions,
see e.g. \cite{Le,Ma}.
For instance, in the inverse spectral problem for 
a one-dimensional Schr\"odinger operator,
 it is known that measuring spectra corresponding to two different boundary
conditions determine the potential uniquely.
Here we consider similar results for the 2-dimensional conductivity
equation assuming that we know measurements on part of the boundary
for two different boundary conditions.

 Let us consider the conductivity equation with the Dirichlet
boundary condition
\beq\label{conduct D}
\nabla\cdotp \sigma\nabla u&=& 0\hbox{ in } \Omega,\\
u|_{\p \Omega}&=&\phi\nonumber
\eeq
and with the Neumann
boundary condition
\beq\label{conduct N}
\nabla\cdotp \sigma\nabla v&=& 0\hbox{ in } \Omega,\\
\nu\cdotp\sigma \nabla v|_{\p \Omega}&=&\psi,\nonumber
\eeq
normalized by $\int_{\p \Om}v\,dS=0$.
Let $\Gamma\subset \p\Om$ be open.
We denote by 
$H_0^{s}(\Gamma)$ the space of functions $f\in H^{s}(\p \Om)$
that are supported on $\Gamma$ and by 
$H^{s}(\Gamma)$ the space of restrictions $f|_\Gamma$ of
 $f\in H^{s}(\p \Om)$.
We define the Dirichlet-to-Neumann map
$\Lambda_\Gamma$ and Neumann-to-Dirichlet map
$\Sigma_\Gamma$ by
\ba
& &\Lambda_\Gamma:H_0^{1/2}(\Gamma)\to H^{-1/2}(\Gamma),
\quad \phi\mapsto (\nu\cdotp\sigma \nabla u)|_{\Gamma},\\
& &\Sigma_\Gamma:H_0^{-1/2}(\Gamma)\to H^{1/2}(\Gamma),
\quad \psi\mapsto v|_{\Gamma}.
\ea
\begin{theorem}\label{Lem: A3} Let $\Gamma\subset \p\Om$ be open.
Then knowing $\p \Om$ and both of the maps
$\Lambda_\Gamma$ and
$\Sigma_\Gamma$ determine the 
equivalence class
\ba
E_{\sigma,\Gamma}=\{\sigma_1\in \Sigma(\Om)& |& 
\hbox{$\sigma_1 = F_*\sigma$, 
$F:\Om\to \Om$ is a $W^{1,2}$-diffeomorphism,}\\
& &F|_{\Gamma}=I\}.
\ea
Moreover, if $\sigma$ is known to be isotropic, it is determined uniquely
by $\Lambda_\Gamma$ and
$\Sigma_\Gamma$.
\end{theorem}

\section{PROOF OF THEOREM \ref{theorem1}} \label{Sec: proofs}

\subsection{Preliminary Considerations}
In the following we identify $\R^2$ and  $\C$  by the 
map $(x^1,x^2)\mapsto x^1+ix^2$
and denote $z=x^1+ix^2$. We use the standard notations
\ba
\p_z=\frac 12(\p_1-i\p_2),\quad \p_{\overline z}=\frac 12(\p_1+i\p_2), 
\ea
where  $\p_j=\p/\p x^j$.  
Below we consider 
$\sigma:\Om\to \R^{2\times 2}$ to be extended as a function 
$\sigma:\C\to \R^{2\times 2}$
by defining $\sigma(z)=I$ for  $z\in \C\setminus \Om$.
In following, we denote $C_0=C_0(\sigma)$.
 For the conductivity $\sigma=\sigma^{jk}$
we define the corresponding Beltrami coefficient (see \cite{Sy1,AP,IM})
\beq\label{mu}
\mu_{1}(z)= 
\frac {-\sigma^{11}(z)+\sigma^{22}(z)-2i\,\sigma^{12}(z)}
 {\sigma^{11}(z)+\sigma^{22}(z)+2\sqrt{\det(\sigma(z))}}.
\eeq
The coefficient $\mu_{1}(z)$ satisfies $|\mu_{1}(z)|\leq \kappa<1$
and is compactly supported.

Next we introduce a $W^{1,2}$-diffeomorphism (not necessarily
 preserving the boundary) 
that transforms the conductivity
to an isotropic one. 

\begin{lemma}\label{lem: 1} 
There is a quasiconformal homeomorphism $F:\C\to \C$
such that
\beq\label{asympt 1}
F(z)=z+\O(\frac 1{z})\quad\hbox{as }|z|\to \infty
\eeq
and such that $F\in W^{1,p}_{loc}(\C;\C)$, $2<p<p(C_0)=\frac {2C_0}{C_0-1}$
for which
\beq\label{isotropic}
(F_*\sigma)(z)=\tilde \sigma(z):=
\det(\sigma(F^{-1}(z)))^{\frac 12}.
\eeq
\end{lemma}

\proof
The proof can be found from   \cite{Sy1} for
$C^3$-smooth conductivities, see also \cite{IM}. 
Because of varying sign conventions, we
 sketch here the proof for readers convenience. We
need to find a quasiconformal map $F$ such that
\beq\label{lassi 1}
DF\,\sigma\,DF^t=\sqrt{\det(\sigma)}\, J_F I
\eeq
where $J_F=\det(DF)$ is the Jacobian of $F$. Denoting by 
$G=[g_{ij}]_{i,j=1}^2$ the inverse of the matrix $\sigma/\sqrt{\det(\sigma)}$
we see that the claim is equivalent to  proving the following:

For any symmetric matrix $G$ with $\det(G)=1$ and $\frac 1 K I
\leq G\leq K I$ there exists a quasiconformal map
$F$ such that
\beq\label{lassi 2}
J_F G= DF^t DF.
\eeq
Next, the non-linear equation
(\ref{lassi 2}) can be replaced in complex notation
by a linear one. Indeed, if $F=u+iv$ then 
(\ref{lassi 2}) is equivalent to
\beq\label{lassi 3}
\nabla v=J G^{-1}\nabla u,\quad \hbox{where }J=
\left(\begin{array}{cc} 
0 & -1      \\ 
1 & 0 
\end{array}\right).
\eeq
This follows readily from the identity
\ba
DF^tJ=\det(DF)\,J\,(DF)^{-1}=JG^{-1}DF^t
\ea
where the latter equality uses (\ref{lassi 2}).
The matrix $J$ corresponds to the multiplication with the imaginary
unit $i$ in complex notation. Denoting by
$R=
\left(\begin{array}{cc} 
1 & 0      \\ 
0 & -1 
\end{array}\right)
$
(the matrix corresponding to complex conjugation) we see that
(\ref{lassi 3}) is equivalent to
\beq\label{lassi 4}
\nabla u+J\nabla v=(G-1)(G+1)^{-1}(\nabla u-J\nabla v).
\eeq
But, $\nabla u+J\nabla v=2\pbz F$ and
$R (\nabla u-J\nabla v)  =2\pz F$
in complex notation and hence (\ref{lassi 4}) becomes 
\beq\label{lassi 5}
& &\p_{\overline z} F=\mu_{1}(z)\p_{z} F
\eeq
where
\ba \mu_{1}=(G-1)(G+1)^{-1}R=(\sqrt {\det \sigma} I-\sigma)
(\sqrt {\det \sigma} I+\sigma )^{-1}R.
\ea
A direct calculation gives
\ba
\mu_{1} =\frac 1{2+g_{11}+g_{22}}\left(\begin{array}{cc} 
g_{11}-g_{22} &  -2g_{12}      \\ 
2g_{12} & g_{11}-g_{22} 
\end{array}\right)
\ea
which shows that the matrix $\mu_{1}=(G-1)(G+1)^{-1}R$ 
corresponds to a multiplication
operator (in complex notation) by the function 
\ba
\mu_{1}(z)=\frac {g_{11}(z)-g_{22}(z)+2ig_{12}(z)}
{2+g_{11}(z)+g_{22}(z)}.  
\ea
This gives (\ref{mu}) since $G^{-1}=\sigma/\sqrt{\det(\sigma)}$.
Since $|\mu_{1}(z)|\leq \kappa <1$ for every $z\in \C$
it is well known by \cite[Thm. V.1, V.2]{Alf} that the equation 
 (\ref{lassi 5}) with asymptotics 
\ba
F(z)=z+{\cal O}(\frac 1z ),\quad\hbox{as }z\to \infty
\ea
has a unique (quasiconformal) solution $F$. The fact that
 $F\in W^{1,p}_{loc}(\C;\C)$, $2<p<\frac {2C_0}{C_0-1}$ follows from
\cite{Astala}.
\Box

In this section we denote by $F=F_\sigma$ the diffeomorphism
determined by Lemma \ref{lem: 1}.  
We also denote $\tilde \Om=F(\Om)$  where
$F$ is as in Lemma \ref{lem: 1}.
Note that (\ref{asympt 1}) implies also that
\beq\label{asympt 2}
F^{-1}(z)=z+\O(\frac 1{|z|})\quad\hbox{as }|z|\to \infty.
\eeq

Later we will use the obvious fact that
the knowledge of map $\Lambda_{\sigma}$ is equivalent
to the knowledge of the Cauchy data pairs
\ba
C_{\sigma}=\{(u|_{\p \Omega},\nu\cdotp \sigma \nabla u|_{\p \Omega})\ |\ u\in H^1(\Om),
\ \nabla\cdotp \sigma \nabla u=0\}.
\ea

In addition to the anisotropic conductivity equation (\ref{conduct})
we consider the corresponding  conductivity equation with isotropic
conductivity. For these considerations,
we observe that if $u$ satisfies equation (\ref{conduct})
and $\widetilde \sigma$ is as in (\ref{isotropic}) then the function
\ba
w(x)=u(F^{-1}(x))\in H^1(\tilde \Om)
\ea
satisfies the isotropic conductivity equation
\beq\label{EQ 1}
& &\nabla\cdotp \widetilde\sigma \nabla w=0\quad \hbox{ in }\tilde \Om,\\
& &w|_{\p \tilde \Omega}=\phi\circ F^{-1}.\nonumber
\eeq
Thus, $\widetilde\sigma$ can be considered as a scalar, isotropic $L^\infty$--smooth
conductivity $\widetilde\sigma I$. We
 continue also the function
$\widetilde\sigma:\tilde \Om\to \R_+$ to a function 
$\widetilde\sigma:\C\to \R_+$
by defining $\widetilde\sigma(x)=1$ for  $x\in \C\setminus \tilde \Om$.

\subsection{Conjugate Functions}
While solving the isotropic inverse problem in \cite{AP}, the interplay of 
the scalar conductivities
$\sigma(x)$ and $\frac 1 {\sigma(x)}$ played a crucial role.
Motivated by this, we define
\ba
\hat \sigma^{jk}(x)=\frac 1{\det(\sigma(x))}\sigma^{jk}(x).
\ea 
Note that for a isotropic conductivity $\hat \sigma=1/\sigma$.

Let now  $F$ be the quasiconformal map defined in  Lemma \ref{lem: 1}
and $\widetilde\sigma=F_*\sigma$ as  in (\ref{isotropic}). We say that
 $\hat w\in H^1(\tilde \Om)$ is a
$\widetilde\sigma$-harmonic conjugate of $w$ if
\beq\label{EQ 2}
& &\p_1\hat w(z)=-\widetilde\sigma(z)\p_2 w(z),\\
& &\p_2\hat w(z)=\widetilde\sigma(z)\p_1 w(z)\nonumber
\eeq
for $z=x^1+ix^2\in \C$.
Using $\hat w$ we define the function $\hat u$ 
that we call the $\sigma$-harmonic conjugate of $u$,
\ba
\hat u(x)=\hat w(F(x)).
\ea
To find the equation governing $\hat u$, it easily follows that
 (c.f. \cite{AP})
 \beq\label{EQ 3}
\nabla\cdotp \frac 1 {\widetilde\sigma} \nabla\hat w=0\quad \hbox{ in }\tilde \Om,
\eeq
and by changing coordinates to $y=F(x)$ we see  that
 $1/\widetilde\sigma=F_*\hat \sigma$.
These facts imply
 \beq\label{EQ 4}
\nabla\cdotp \hat \sigma \nabla\hat u=0\quad \hbox{ in }\Om.
\eeq
Thus $\hat u$ is the $\widehat\sigma$-harmonic conjugate
function of $u$ and we have 
\beq\label{EQ 5}
\nabla \hat u=J \sigma \nabla u,\quad
\nabla  u=J\widehat \sigma \nabla \hat u.
\eeq
 Since $u$ is a solution of the 
conductivity equation if and only if $u+c$, $c\in \C$, is 
solution, we see from (\ref{EQ 5}) that the Cauchy data 
pairs $C_{\sigma}$ determine the pairs  $C_{\hat \sigma}$ 
and vice versa. Thus we get, almost free, that $\Lambda_\sigma$
determines  $\Lambda_{\hat \sigma}$, too.

Let us next consider the function
\beq\label{EQ 6}
f(z)=w(z)+i\hat w(z).
\eeq
By \cite{AP}, it satisfies the pseudo-analytic equation
of second type,
\beq\label{EQ 7}
\pbz f=\tilde \mu_{2}\, \overline {\pz f}
\eeq
where 
\beq\label{EQ 7b}
\tilde \mu_{2}(z)=\frac {1-\widetilde\sigma(z)}{1+\widetilde\sigma(z)},\quad 
|\tilde \mu_{2}(z)|\leq 
\frac {C_0-1}{C_0+1}<1.
\eeq
Using this Beltrami coefficient,
we define $\mu_2=\tilde \mu_2\circ F$.

We will need the following:

\begin{lemma}\label{lemma21}
Let $g=f\circ F$ where 
$F:\Om\to \widetilde \Om$ is a quasiconformal homeomorphism and $f$
is a quasiregular map satisfying
\beq\label{L99}
\pbz f=\tilde \mu_2 \overline{\pz f}\quad\
\hbox{and}\quad
\pbz F=\mu_1 \pz F,
\eeq
where $\tilde \mu_2=\mu_2\circ F^{-1}$ and $\mu_1$ satisfies
$|\mu_j|\leq \kappa<1$ and $\mu_2$ is real.
Then $g$ is quasiregular and satisfies
\beq\label{L100}
\pbz g= \nu_1{\pz g}+ \nu_2 \overline{\pz g},
\eeq
 where
\beq\label{L101}
\nu_1=\mu_1 \frac {1-\mu_2^2}{1-|\mu_1|^2\mu_2^2}
,\quad\
\hbox{and}\quad
\nu_2=\mu_2 \frac {1-|\mu_1|^2}{1-|\mu_1|^2\mu_2^2}.
\eeq
Conversely, if $g$ satisfies (\ref{L100}) with
$\nu_2$ real and $|\nu_1|+|\nu_2|\leq \kappa' <1$
then there exists unique $\mu_1$ and $\mu_2$
such that (\ref{L101}) holds and
 $f=g\circ F^{-1}$ satisfies (\ref{L99}).
\end{lemma}

\proof
We apply the chain rule
\ba
\p(f\circ F)&=&
(\p f)\circ F\cdotp\p F+(\bar \p f)\circ F\cdotp\bar{\bar \p  F},\\
\bar \p(f\circ F)
&=&
(\p f)\circ F\cdotp\bar \p F+(\bar \p f)\circ F\cdotp\bar{ \p  F},
\ea
and obtain
\ba
 \nu_1\pz g+ \nu_2 \overline{\pz g}=
\p f\circ F\cdotp \p F \cdotp (\nu_1+\nu_2\mu_1\mu_2)+
\overline {\p f\circ F}\cdotp {\overline\p F}\cdotp (\nu_2+\nu_1\overline \mu_1\mu_2)
\ea
and
\ba
\pbz g= \mu_1 \cdotp \p f\circ F\cdotp \p F 
+
\mu_2 \cdotp \overline {\p f\circ F}\cdotp {\overline\p F}. 
\ea
Hence, if $\mu_1,\mu_2,\nu_1$, and $\nu_2$ are related
so that
\ba
\mu_1=\nu_1+ \nu_1\mu_2,\quad
\mu_2=\nu_2+\overline \nu_1\mu_2,
\ea
we see that (\ref{L100}) and (\ref{L101}) are satisfied.

It is not difficult to see that for each $\nu_1$ and
$\nu_2$ (\ref{L101}) has a unique solution
$\mu_1,\mu_2$ with $|\mu_j|\leq \kappa'<1$, $j=1,2$.
Again, the general theory of quasiregular maps \cite{Alf} implies
that (\ref{L99}) has a solution and the factorization
$g=f\circ F$ holds.
\Box

Note that (\ref{L101}) implies that 
\beq\label{L103}
|\nu_1|+|\nu_2|= \frac  {|\mu_1|+|\mu_2|}{1+|\mu_1|\,|\mu_2|}
\leq 
\frac {2\kappa}{1+\kappa^2}<1.
\eeq

Lemma \ref{lemma21} has the following important corollary,
that is the main goal of this subsection.

\begin{corollary}\label{cor21}
If $u\in H^1(\Om)$ is a real solution
of the conductivity equation (\ref{conduct}),
there exists  $\hat u\in H^1(\Om)$, unique up to
a constant, such that $g=u+i\hat u$ satisfies
(\ref{L100}) where
\beq\label{L104}
\nu_1= \frac {\sigma^{22}-\sigma^{11}-2i\sigma^{12}}
{1+\hbox{\rm tr}\, \sigma+\det(\sigma)}
,\quad\
\hbox{and}\quad
\nu_2= \frac {1-\det(\sigma)}
{1+\hbox{\rm tr}\, \sigma+\det(\sigma)}.
\eeq
Conversely, if 
$\nu_1$ and $\nu_2$, $|\nu_1|+|\nu_2|\leq \kappa' <1$ 
 are as in Lemma \ref{lemma21}
then there are unique $\sigma$ and $\hat\sigma$ such
that for any solution
 $g$ of (\ref{L100})
 $u=\re g$ and $\hat u=\im g$ satisfy the
conductivity equations
\beq\label{cond eq 2}
\nabla\cdotp \sigma\nabla u=0,\quad\hbox{and}\quad
\nabla\cdotp \hat \sigma\nabla \hat u=0.
\eeq
\end{corollary}

\proof Since $g=f\circ F$ where $
f=w+i\hat w$ according to (\ref{EQ 6}),
we obtain immediately the existence of $\hat u=\hat w\circ F$.
Thus we need only to calculate $\nu_1$ and $\nu_2$
in terms of $\sigma$. Note that by 
(\ref{mu}), 
\beq\label{L106}
|\mu_1|^2= \frac {\text{tr}\,(\sigma)-2\det(\sigma)^{1/2}}
{\text{tr}\,(\sigma)+2\det(\sigma)^{1/2}}.
\eeq
We recall that
\beq\label{L107}
\mu_2=
\frac {1-\det(\sigma)^{1/2}}
{1+\det(\sigma)^{1/2}}
\eeq
and thus
\ba
1-|\mu_1|^2\mu_2^2=
\frac {4(\det(\sigma)^{1/2}\text{tr}\,(\sigma)+
(1+\det(\sigma))\det(\sigma)^{1/2})}
{(1+\det(\sigma)^{1/2})^2(\text{tr}\,(\sigma)+
2\det(\sigma)^{1/2})}
\ea
which readily yields (\ref{L104}) from (\ref{L101}).

Note that since $\nu_1$ and $\nu_2$ uniquely
determine $\mu_1$ and $\mu_2$, they by (\ref{L106})
and  (\ref{L107}) also determine $\det (\sigma)$
and tr$(\sigma)$. After observing this, it is clear from
(\ref{L104}) that $\sigma$ is uniquely determined by $\nu_1
$ and $\nu_2$.
\Box

Now one can write equations (\ref{EQ 5}) in
 more explicit
 form
\beq\label{EQ 9}
\tau\cdotp \nabla \hat u|_{\p \Om}=\Lambda_\sigma (u|_{\p \Om})
\eeq
where $\tau=(-\nu_2,\nu_1)$ is a unit tangent vector
of $\p \Om$. As $\re g|_{\p \Om}=u|_{\p \Om}$
and  $\im g|_{\p \Om}=\hat u|_{\p \Om}$, we see that
 $\Lambda_\sigma$ determines
the $\sigma$-Hilbert transform $\H_\sigma$ defined by
\beq\label{EQ 9b}
\H_\sigma&:& H^{1/2}(\p \Om)\to H^{1/2}(\p \Om)/\C,\\
\nonumber
& &\re g|_{\p \Om}\mapsto \im g|_{\p \Om}+\C.
\eeq
Put yet in another terms, for $u,\hat u\in H^{1/2}(\p \Om)$,
$\hat u={\cal H}_\sigma u$ if and only if the map
$g(\xi)=(u+i\hat u)(\xi),$ $\xi \in \p \Om$, extends
to $\Om$ so that (\ref{L100}) is satisfied.

Summarizing the previous results, we have 
\begin{lemma}\label{lem: 2} 
The Dirichlet-to-Neumann map $\Lambda_\sigma$
determines the maps $\Lambda_{\hat \sigma}$ and $\H_\sigma$.
\end{lemma}

\subsection{Solutions of Complex Geometrical Optics}

Next we consider exponentially growing solutions,
i.e., solutions of complex geometrical optics originated
by Calder\'on for linearized inverse problems and by Sylvester and
Uhlmann for non-linear inverse problems.
In our case, we seek solutions
$G(z,k)$, $z\in \C\setminus \Om$, $k\in \C$ satisfying
\beq\label{EQ 10}
& &\pbz G(z,k)=0 \quad \hbox{for }z\in \C\setminus \overline \Om,\\
\label{EQ 10 asym}
& &G(z,k)=e^{ikz}(1+\O_k(\frac 1{z})),\\
\label{EQ 10 bnd}
& &\im G(z,k)|_{z\in \p \Om}=\H_\sigma
(\re G(z,k)|_{z\in \p \Om}).
\eeq
Here $\O_k(h(z))$ means a function of $(z,k)$ that satisfies
$|\O_k(h(z))|\leq C(k)|h(z)|$ for all $z$ with some constant $C(k)$
depending only on
$k\in \C$.
For the conductivity $\widetilde\sigma$ we consider the 
corresponding exponentially growing solutions
 $W(z,k)$, $z\in \C$, $k\in \C$ where
 \beq\label{EQ 11}
& &\pbz W(z,k)=\tilde \mu_{2}(z)\overline {\pz W(z,k)}, \quad \hbox{for }z\in \C,\\
\label{EQ 11 b}
& &W(z,k)=e^{ikz}(1+\O_k(\frac 1 z)).
\eeq
Note that in this stage, $z\mapsto G(z,k)$ is defined only
in the exterior domain $\C\setminus \Om$ but $z\mapsto W(z,k)$ in the whole complex plane.
These two solutions are closely related:

\begin{lemma}\label{lem: 3} For all  $k\in \C$ we have:

i. The system (\ref{EQ 11}) has a unique solution
$W(z,k)$ in $ \C$.\\

ii. The system (\ref{EQ 10}--\ref{EQ 10 bnd}) has a unique solution
$G(z,k)$ in $\C\setminus \Om$.\\

iii. For $z\in \C\setminus \Om$  we have
\beq\label{EQ 11b}
G(z,k)=W(F(z),k).
\eeq
\end{lemma}

\proof
For the claim i. we refer to \cite[Theorem 4.2]{AP}.

Next we consider ii. and iii. simultaneously.
Assume that $G(z,k)$ is a solution
of (\ref{EQ 10}--\ref{EQ 10 bnd}). By Lemma \ref{lemma21} and 
boundary condition (\ref{EQ 10 bnd}) we 
see that equation
\ba
& &\bar \p h(z,k)
= \nu_1\pz h+\nu_2 \overline{\pz h}, \quad \hbox{in }\Om, \\
& &h(\cdotp,k)|_{\p \Omega}=G(\cdotp,k)|_{\p \Omega}
\ea
has a unique solution where $\nu_1$ and $\nu_2$ are given 
in (\ref{L104}). 

 Let
\beq\label{extensions}
H(z,k)=\begin{cases}
G(z,k)&\quad \hbox{for } z\in \C\setminus \Om\\
h(z,k)&\quad \hbox{for } z\in \Om\end{cases}
\eeq
and $\tilde H(z,k)=H(F^{-1}(z),k)$. Then 
$\tilde H(z,k)$ satisfies equations 
\ba
& &\pbz \tilde H(z,k)=0, \quad \hbox{for }z\in \C\setminus \tilde \Om,\\
& &\pbz \tilde H(z,k)=\tilde \mu_{2}(z)\overline {\pz \tilde H(z,k)}, \quad \hbox{for }z\in 
\tilde \Om,
\ea
and traces from both sides of $\p\tilde \Om$ coincide.
Thus $\tilde H(z,k)$ satisfies  equation
in (\ref{EQ 11}).
Now (\ref{asympt 2}) and (\ref{EQ 10 asym}) yield
that 
\beq\label{EQ 11c}
\tilde H(z,k)&=&H(F^{-1}(z),k)\\
&=& \nonumber
\exp(ikF^{-1}(z))(1+\O_k(\frac 1{1+|F^{-1}(z)|}))\\
&=& \nonumber
\exp(ikz)(1+\O_k(\frac 1{1+|z|}))
\eeq
showing that $\tilde H$ satisfies
(\ref{EQ 11}--\ref{EQ 11 b}). Thus
 by i., $\tilde H(z,k)=W(z,k)$.
This proves both ii. and iii.
\Box

\subsection{Proof of Theorem \ref{theorem1}}
As $G(z,k)$ is the unique solution of (\ref{EQ 10}--\ref{EQ 10 bnd})
and the  operator appearing in boundary condition 
  (\ref{EQ 10 bnd}) is  known,
Lemmata \ref{lem: 2} and  \ref{lem: 3}
imply the following:

\begin{lemma}\label{lem: 3b} The Dirichlet-to-Neumann map
$\Lambda_\sigma$ determines
$G(z,k)$, $z\in \C\setminus \Om$, $k\in \C$.
\end{lemma}

Next we use this to find the diffeomorphism $F_\sigma$ outside $\Om$.

\begin{lemma}\label{lem: 4} The Dirichlet-to-Neumann map
$\Lambda_\sigma$ determines the values the restriction
$F_\sigma|_{\C\setminus \Om}$.
\end{lemma}

\proof
By (\ref{EQ 11b}), $G(z,k)=W(F(z),k)$,
where $W(z,k)$ is the exponentially growing solution corresponding
to the isotropic conductivity $\widetilde\sigma$. Thus by applying
the sub-exponential growth results for such solutions,
\cite[Lemma 7.1 and Thm. 7.2]{AP}, we have representation
\beq\label{EQ 17}
W(z,k)=\exp(ik\varphi(z,k))
\eeq
where
\beq\label{EQ 18}
\lim_{k\to \infty} \sup_{z\in \C}|\varphi(z,k)-z|=0.
\eeq
As $F(z)=z+{\cal O}(1/z)$, and  $G(z,k)=W(F(z),k)$ we have
\beq\label{EQ 19}
\lim_{k\to \infty} \frac{\log G(z,k)}{ik}=
\lim_{k\to \infty} \varphi(F(z),k)=F(z).
\eeq
By Lemma \ref{lem: 3b} we know the values
of limit (\ref{EQ 19}) for any $z\in \C\setminus\Omega$. Thus the claim is
proven. \Box 

We are ready to prove Theorem \ref{theorem1}.

\proof
As we know $F|_{\C\setminus \Om}\in W^{1,p}$, $2<p<p(C_0)$,
 we in particular know
$\tilde \Om=\C\setminus (F(\C\setminus \Om)).$ 
When $u$ is the solution of conductivity equation (\ref{conduct})
with Dirichlet boundary value $\phi$ and $w$  is the solution
of (\ref{EQ 1}) with Dirichlet boundary value $\widetilde \phi
=\phi\circ h$, where $h=F^{-1}|_{\p \tilde \Om}$ we see that 
\beq\label{kaava 1 apu}
\int_{\p \widetilde \Om}\widetilde \phi  
\Lambda_{\widetilde\sigma}\widetilde \phi\, dS =
Q_{\widetilde\sigma,\widetilde\Om}(w)=
Q_{\sigma,\Om}(u)=\int_{ \p \Om}\phi  
\Lambda_{\sigma}\phi\, dS. 
\eeq
Here, the second identity is justified by the fact that $F$
is quasiconformal and hence satisfies (\ref{Kari 9}). 
Since $
\Lambda_{\sigma}$ and 
$\Lambda_{\widetilde\sigma}$ are symmetric, this implies
\beq\label{kaava 1}
\int_{\p \widetilde \Om}\widetilde \psi  
\Lambda_{\widetilde\sigma}\widetilde \phi\, dS =
\int_{ \p \Om}\psi  
\Lambda_{\sigma}\phi\, dS 
\eeq
for any $\widetilde \psi ,\widetilde \phi\in H^{1/2}(\p
\widetilde \Om)$ and  $ \psi , \phi\in H^{1/2}(\p
 \Om)$ are related by $\widetilde \psi= \psi\circ h$
and $\widetilde \psi= \psi\circ h$. Note that
 $\phi\in H^{1/2}(\p \Om)$ if and only if
$\widetilde \phi=\phi \circ h\in H^{1/2}(\p \widetilde\Om)$.
To see this, extend $\phi$ to a $H^1(\Om)$ function
and after that define $\widetilde \phi$ in the interior
of $\widetilde\Om$ by $\widetilde \phi=\phi\circ F^{-1}$.
Now
\ba 
||\nabla \phi||_{L^2(\Om)}^2
\sim 
\int_{\Om } \nabla \phi\cdotp \sigma \overline {\nabla \phi}\,dx
\sim 
\int_{\widetilde\Om } \nabla \widetilde\phi\cdotp 
\widetilde\sigma \overline {\nabla \widetilde\phi}\,dx
\sim ||\nabla \widetilde\phi||_{L^2(\widetilde\Om)}^2
\ea
and hence
\ba 
||\phi||_{H^{1/2}(\p \Om)}^2
\sim 
|| \phi||_{H^1(\Om)}^2
\sim 
||\widetilde \phi||_{H^1(\widetilde \Om)}^2
\sim 
||\widetilde \phi||_{H^{1/2}(\p \widetilde \Om)}^2.
\ea

As we know $F|_{\C\setminus \Om}$ and $\Lambda_\sigma$,
we can find $\Lambda_{\widetilde \sigma}$ using formula
(\ref{kaava 1}).
By \cite{AP}, the map $\Lambda_{\widetilde \sigma}$ 
determines uniquely the conductivity $\widetilde\sigma$
on $\tilde \Om$ in a constructive manner. 

Knowing $\Om$, $\tilde \Om$, and the boundary value
$f=F|_{\p \Om}$ of the map $F:\Om\to \tilde \Om$, we next 
construct a sufficiently smooth diffeomorphism
$H:\tilde \Om\to \Om$. First, by the Riemann mapping
theorem we can map $\Om$ and $\tilde \Om$ to the unit
disc $\Disc$ by the conformal maps $R$ and $\tilde R$, respectively.
Now 
\ba
G=\tilde R\circ F\circ R^{-1}:\Disc\to\Disc
\ea  
is a quasiconformal map and since we know $R$ and $\tilde R$,
we know the function $g=G|_{\p \Disc}$ mapping
$\p \Disc$ onto itself. The map $g$ is quasisymmetric 
(cf. \cite{Alf}) and by 
Ahlfors-Beurling extension theorem \cite[Thm. IV.2]{Alf} 
it has a quasiconformal
extension ${\mathcal AB}(g)$ mapping $\overline \Disc$
onto itself. Note that one can obtain  ${\mathcal AB}(g)$
from $g$ constructively by an explicit formula 
 \cite[p. 69]{Alf}. Thus we may find a quasiconformal diffeomorphism
$H=R^{-1} \circ[{\mathcal AB}(g)]^{-1}\circ \tilde R$, $H:\tilde \Om\to \Om$ 
satisfying $H|_{\p \tilde \Om}=
F^{-1}|_{\p \tilde \Om}$. 

Combining the above results, we can find
$H_*\widetilde\sigma$ that is a representative of the equivalence class
$E_\sigma$. 
\Box

 In the above proof the Riemann mappings can not be found as explicitly
as the Ahlfors-Beurling extension. However, there are
numerical packages for approximative construction of
Riemann mappings, see e.g. \cite{webbi}.

\section{PROOFS OF CONSEQUENCES OF MAIN RESULT}
\label{sec: proof of con}

Here we give proofs of Theorems \ref{Lem: A1}--\ref{Lem: A3}.

{\bf Proof of  Theorem \ref{Lem: A1}.}
Let $F:\R_-^2=\R+i\R_-\to \Disc $ be the M\"obius transform
\ba
F(z)=\frac {z+i}{z-i}.
\ea
Since this map is conformal, we see that
$C_0(F_*\sigma)=C_0(\sigma).$
Let $\tilde \sigma=F_*\sigma$ be the conductivity in $\Disc $.
Then $\Lambda_{\tilde \sigma}\phi$ is determined as in
(\ref{kaava 1}) for all $\phi\in C^\infty_0(\p \Disc \setminus \{1\})$.

Since  $\Lambda_{\tilde \sigma}1=0$ and functions
$\C\oplus C^\infty_0(\p \Disc \setminus \{1\})$ are dense
in the space $H^{1/2}(\p \Disc )$, we see that 
$\Lambda_{\sigma}$ determines the Dirichlet-to-Neumann map
$\Lambda_{\tilde \sigma}$ on $\p \Disc $. Thus we 
can find the equivalence class of the conductivity on $\Disc $.
Pushing these conductivities forward with $F^{-1}$ to $\R^2_-$, we obtain the claim.
\Box

{\bf Proof of  Theorem \ref{Lem: A2}.}
Let $F:S\to \Disc \setminus \{0\}$ be the conformal map
such that
\ba
\lim_{z\to \infty}F(z)=0.
\ea
Again, since this map is conformal we have for $\tilde \sigma=F_*\sigma$
the equality
$C_0(\sigma)=C_0(F_*\sigma).$
Moreover, if $u$ is a solution of (\ref{conduct S}),
we have that $w=u\circ F^{-1}$ is solution of
\beq\label{conduct S2}
\nabla\cdotp \tilde \sigma\nabla w&=& 0\hbox{ in } \Disc \setminus \{0\},
\\
w|_{\p \Disc}&=&\phi\circ F^{-1},\nonumber\\
w&\in& L^\infty( \Disc)\nonumber.
\eeq
Since set $\{0\}$ has capacitance zero in $\Disc $, we 
see that $w=W|_{\Disc \setminus \{0\}}$ where
\beq\label{conduct S3}
\nabla\cdotp \tilde \sigma\nabla W&=& 0\hbox{ in } \Disc ,
\\
W|_{\p  \Disc}&=&\phi\circ F^{-1}.\nonumber
\eeq
Since $F$ can be constructed via the Riemann mapping theorem,
we see that $\Lambda_\sigma$ determines
$\Lambda_{\tilde \sigma}$ on $\p \Disc $ and thus
the equivalence class $E_{\tilde \sigma}$. When
$\tilde F:\Disc \to \Disc $ is a boundary preserving 
diffeomorphism, we see that $F^{-1}\circ\tilde F\circ F$
defines a diffeomorphism $\overline S\to \overline S$. 
Since we have determined
the conductivity $\tilde \sigma$ up to a 
 boundary preserving 
diffeomorphism, the claim follows easily.
\Box

{\bf Proof of  Theorem \ref{Lem: A3}.}
Let $\Disc\subset \C$ be the unit disc and
$\Disc_+=\{z\in \Disc\ | \ \re z>0\}$.
Let $F:\Om\to \Disc_+$ be a Riemann mapping 
such that
\ba
\Disc_+\subset \R\times \R_+,\quad F(\Gamma)=\p \Disc_+\setminus
(\R\times \{0\}),\quad F(\p \Om \setminus \Gamma)=\p \Disc_+\cap
(\R\times \{0\}).
\ea
Let $\eta:(x^1,x^2)\mapsto (x^1,-x^2)$ 
and define $\Disc_-=\eta (\Disc_+)$, 
and $\tilde \sigma=F_*\sigma$. Let
\ba
\hat \sigma(x)=
\begin{cases}
\sigma(x) &\quad \hbox{for } x\in \Disc_+,
\\
(\eta_*\sigma)(x) &\quad \hbox{for } x\in \Disc_-.\end{cases}
\ea
Consider equation
\beq\label{conduct D+-}
\nabla\cdotp \hat \sigma\nabla w&=& 0\hbox{ in } \Disc.
\eeq
Using formula (\ref{kaava 1}) we see that
$F$ and $\Lambda_\Gamma$ determine the corresponding
map $\Lambda_{F(\Gamma)}$ for $\hat \sigma$.
Similarly, we can find  $\Sigma_{F(\Gamma)}$ for $\hat \sigma$.

Then $\Lambda_{F(\Gamma)}$ determines the Cauchy data
on the boundary for the solutions 
of (\ref{conduct D+-}) for which
$w\in H^1(\Disc)$, $w=-w\circ \eta$. On the other hand, 
 $\Sigma_{F(\Gamma)}$ determines the Cauchy data
on the boundary of the solutions 
of (\ref{conduct D+-}) for which
$w\in H^1(\Disc)$ and $w=w\circ \eta$.
Now each solution $w$ of (\ref{conduct D+-}) can be
written as a linear combination 
\ba
w(x)=\frac 12(w(x)+w(\eta(x)))+\frac 12(w(x)-w(\eta(x))).
\ea
Thus the maps  $\Lambda_{F(\Gamma)}$  and $\Sigma_{F(\Gamma)}$
together  determine $C_{\hat \sigma}$, and hence we
can find $\hat \sigma$ up to a diffeomorphism.
We can choose
a representative $\hat \sigma_0$ of the equivalence class 
$E_{\hat \sigma}$ such that
 $\hat \sigma_0= \hat \sigma_0\circ \eta$. 
In fact, choosing a symmetric
Ahlfors-Beurling extension in the construction
given in the proof of Theorem \ref{theorem1},
we obtain such a conductivity.
Pushing the conductivity  $\hat \sigma_0$ 
from  $\Disc_+$ to $\Om$ with $F^{-1}$, we 
obtain the claim.
\Box

\bibliographystyle{amsalpha}

\end{document}